\documentclass[12pt]{amsart}
\allowdisplaybreaks[4]
\usepackage{amssymb,color}
\usepackage[square,compress,comma, numbers,sort]{natbib}

\usepackage[colorlinks=true, citecolor=blue, linkcolor=blue]{hyperref}
\usepackage{amsthm}
\usepackage{mathtools}

\definecolor{c20}{rgb}{0.,0.7,0.}
\definecolor{c30}{rgb}{0.,0.,1.}
\definecolor{c40}{rgb}{1,0.1,0.7}
\definecolor{c50}{rgb}{1,0,0}
\definecolor{c60}{rgb}{0,0.9,0.1}

\newcommand{\E}[1]{\mathbb{E} \left(#1\right)}

\newcommand{\pk}[1]{\mathbb{P} \left(#1 \right) }

\newcommand{\R}{\mathbb{R}}

\newcommand{\BQN}{\begin{eqnarray}}
\newcommand{\EQN}{\end{eqnarray}}
\newcommand{\BQNY}{\begin{eqnarray*}}
\newcommand{\EQNY}{\end{eqnarray*}}

\def\eL#1{\textcolor{c30}{#1}}
\def\eL#1{{#1}}

\def\K1#1{\textcolor{cyan}{#1}}
\def\K1#1{#1}

\newcommand{\BS}{\begin{sat}}
\newcommand{\ES}{\end{sat}}
\newcommand{\BT}{\begin{theo}}
\newcommand{\ET}{\end{theo}}
\newcommand{\BK}{\begin{korr}}
\newcommand{\EK}{\end{korr}}

\newcommand{\BD}{\begin{de}}
\newcommand{\ED}{\end{de}}

\newcommand{\BIT}{\begin{itemize}}
\newcommand{\EIT}{\end{itemize}}
\newcommand{\BDI}{\begin{description}}
\newcommand{\EDI}{\end{description}}

\newcommand{\BRM}{\begin{remarks}}
\newcommand{\ERM}{\end{remarks}}

\newcommand{\BEL}{\begin{lem}}
\newcommand{\EEL}{\end{lem}}

\newtheorem{theo}{Theorem}[section]
\newtheorem{sat}[theo]{Proposition}
\newtheorem{de}[theo]{Definition}
\newtheorem{lem}[theo]{Lemma}

\newtheorem{korr}[theo]{Corollary}

\newtheorem{remarks}[theo]{Remarks}

\newcommand{\prooftheo}[1]{ \textbf{Proof of Theorem} \ref{#1} }

\newcommand{\COM}[1]{}

\newcommand{\QED}{\hfill $\Box$ \\}

\def\rw{\rightarrow}

\def\IF{\infty}

\newcommand{\expon}[1]{\exp\left( #1 \right)}

\def\eHH#1{\textcolor{c50}{#1}}

\def\eHH#1{#1}

\def\vf{\sigma^2}

\begin{document}

\title
{  The time of ultimate recovery in Gaussian risk model}

\author{Krzysztof D\c{e}bicki}
\address{Krzysztof D\c{e}bicki, Mathematical Institute, University of Wroc\l aw, pl. Grunwaldzki 2/4, 50-384 Wroc\l aw, Poland}
\email{Krzysztof.Debicki@math.uni.wroc.pl}
\author{Peng Liu}
\address{Peng Liu, Department of Actuarial Science, Faculty of Business and Economics, University of Lausanne, UNIL-Dorigny 1015 Lausanne, Switzerland}
\email{peng.liu@unil.ch}

\bigskip

 \maketitle
\bigskip
{\bf Abstract:} We analyze the distance $\mathcal{R}_T(u)$
between the first and the last passage time
of $\{X(t)-ct:t\in [0,T]\}$ at level $u$ in time horizon $T\in(0,\infty]$,
where $X$ is a centered Gaussian process with stationary increments
and $c\in\R$, given that the first passage time occurred before $T$.
Under some tractable assumptions on $X$, we
find $\Delta(u)$ and $G(x)$ such that
$$\lim_{u\to\infty}\pk{\mathcal{R}_T(u)>\Delta(u)x}=G(x),$$
for $x\ge 0$.
We distinguish
two scenarios: $T<\infty$ and $T=\infty$, that lead to
qualitatively different asymptotics.
The obtained results
provide
exact asymptotics
of the ultimate recovery time after the ruin in Gaussian risk model.



\section{Introduction}

For given threshold $u>0$ and time horizon $T\in(0,\infty]$, let
$$\tau_T(u):=\inf\left\{t\geq 0: X(t)-ct>u, 0\leq t\leq T\right\}$$
and
$$\mathcal{T}_T(u):=\sup\left\{t\geq 0: X(s)-cs>u,  0\leq t\leq T \right\},$$
be {\it the first} and {\it the last passage time} of process $X(t)-ct,t\ge0$ at level $u$ respectively,
with convention that $\inf\emptyset=\IF$ and $\sup\emptyset=0$.

The analysis of properties of $\tau_T(u)$
and $\mathcal{T}_T(u)$, due to their obvious importance in extreme value theory
of stochastic processes, attracted substantial interest, being
additionally stimulated by relations of passage times with
important problems in applied probability.
More specifically, in risk theory
$\tau_T(u)$
and $\mathcal{T}_T(u)$ have the interpretation
as {\it the first} and {\it the last ruin time}
of the risk reserve process
$S(t):=u+ct-X(t)$,
where $u>0$ is the initial capital, $c>0$ is the premium rate and $X(t),t\ge0$ is the
accumulated claim amount in interval $(0,t]$; see e.g. \cite{ChS05,HJ13}.

In this contribution we suppose that $X$ is
a centered Gaussian process with stationary increments, a.s.
continuous sample paths, $X(0)=0$ and $c>0$.
We note that in the context of risk theory,
there are strong application-based and  theoretical reasons for modelling accumulated claim amount by
Gaussian processes with stationary increments.
On one hand the family of Gaussian processes
provides flexibility in the adjustment of suitable correlation model, since it covers
wide range of correlation structures.
On the other hand, there are theoretical results that legitimate approximation of
the accumulated claim amount in highly aggregated models by Gaussian processes;
see the celebrated work by Iglehart \cite{Igl69} for the Brownian approximation and, e.g.,
\cite{Mic98, HJ13, KEP2015} for more general models, including e.g. fractional Brownian motion approximations.

This contribution is devoted to analysis of the
distribution of \[
\mathcal{R}_T(u):=\mathcal{T}_T^*(u)-\tau_{T}^*(u),
\]
where
\BQN\label{tau}
\left(\tau_T^*(u),\mathcal{T}_T^*(u)\right):= \left(\tau_T(u),\mathcal{T}_T(u)\right)\Bigl\lvert(\tau_T(u)<\IF),
\EQN
both for $T\in(0,\infty)$ and $T=\infty$.
Referring again to risk theory,
$\mathcal{R}_T(u)$
has the interpretation as {\it the ultimate time to recovery}, which is the
difference between the last and first ruin time, under the condition that ruin occurred; see also
\cite{Fro04,Li08} and references therein.
We note that $\mathcal{R}_T(u)$ is
also closely related to the so-called {\it Parisian ruin time},
which is the first time that the length of the consecutive excursion period of the surplus process $S$ under level 0
exceeds a pre-specified time threshold; see, e.g., \cite{Che97}, \cite{Pal13}, \cite{DHJ15} and \cite{DHJ16},
with straightforward observation that
$\mathcal{R}_T(u)$ gives an upper bound
for the appropriately chosen pre-specified time period in Parisian model.
Another related notion is {\it the cumulative ruin time}, which is based on the total time spent below 0 ({\it in red})
by the underlying risk process; see, e.g., \cite{Gu17}.
Clearly, $\mathcal{R}_T(u)$ is greater than the corresponding occupation time.

For  $T=\IF$, the asymptotics of the distribution of conditional first and last ruin times
in Gaussian risk context were studied in, e.g., \cite{Zhang082} and \cite{Zhang081}; see also
\cite{HJ13}, \cite{PengLith} and \cite{KEP2015} for related
$\gamma-$reflected Gaussian models.
Specifically, under some tractable assumptions on $X$,
the following asymptotics was found in  \cite{Zhang082,Zhang081,KEP2015}:
\BQN\label{joint}
\left(\frac{\tau_\IF^*(u)-ut_u}{\sigma(u\eL{t_u})},\frac{\mathcal{T}_\IF^*(u)-ut_u}{\sigma(u\eL{t_u})}
\right)\stackrel{d}{\rw} \left({\rm C}\mathcal{N},{\rm C} \mathcal{N}\right), \quad u\rw\IF,
\EQN
where $\sigma^2(t)=Var(X(t))$, $t_u=\arg\inf_{t>0}\frac{u(1+ct)}{\sigma(ut)}$,
$\mathcal{N}\sim N(0,1)$ and ${\rm C}$ is some known constant.
However, the above result is too crude in order to deal with $\mathcal{R}_\infty(u)$, since
it follows straightforwardly from (\ref{joint}), that
$$\frac{\mathcal{T}_\IF^*(u)-\tau_\IF^*(u)}{\sigma(ut_u)}\rw 0, \quad u\rw \IF$$
in probability; see also Corollary 4 in \cite{Zhang081}.
This implies that $\sigma(ut_u)$ acting as a scaling function is too big, so
in order to get a nontrivial result one has to
scale by asymptotically smaller function.

The main results of this contribution provide $\Delta(u)$
and distribution function $G(\cdot)$ such that
$$\lim_{u\to\infty}\pk{\mathcal{R}_T(u)>x\Delta(u)}=G(x),$$
for $x\ge 0$, both for $T\in(0,\infty)$ and $T=\infty$.\\
As it is shown in Theorems \ref{TH4}, \ref{TH5}, both $\Delta(u)$
and $G(x)$
depend  on $T$ and on the
local behavior of variance function of $X$, which leads to several
scenarios. Interestingly, the limit function $G$ is given in terms of
{\it generalized Pickands-Piterebarg-type} constants.
In order to obtain the main results of this contribution we
accommodate to our needs recently developed {\it uniform double-sum method}
applied for relevant continuous functionals; see \cite{KEP2016}.

\COM{
In this paper, we focus on the investigation of  asymptotic distribution of  $\frac{\mathcal{T}_T^*(u)-\tau_T^*(u)}{A(u)}$
with appropriately chosen scaling function $A(u)$ as $u\rw\IF$ for $T\in (0,\IF]$.
Our main results below show that the appropriate scaling functions are case-specific and the asymptotic
distribution functions are given Pickands-Piterbarg constants, which are new constants that first appear.
This kind of distributions makes the results much more complicated than normal distribution  obtained in
\cite{Zhang082} and \cite{Zhang081}. The technique used in this paper is the uniform double-sum method in
\cite{KEP2016} which is the development of the celebrated double-sum method in \cite{pit96}.
However, in order to derive our results, we partially extend the uniform approximation for general functional
of Gaussian random fields in \cite{KEP2016}.\\}


Organization of the paper:  Section 2 is devoted to introduction of
notation and presentation of main results. In section 3 we present  proofs of the main results.
\section{Main results}
In this section we provide main results of this contribution, which is the limit theorem for
\[
\mathcal{R}_T(u)=\mathcal{T}_T^*(u)-\tau_{T}^*(u),
\]
as $u\to\infty$, for $T\in(0,\infty]$.
Due to specific asymptotic nature of $\mathcal{R}_T(u)$
we distinguish two separate scenarios: infinite-time horizon ($T=\infty$) and
finite-time horizon ($T<\infty$).

\subsection{Infinite-time horizon} Suppose that $T=\infty$.
Consider a centered Gaussian process $X$ with continuous trajectories, stationary increments and
variance function $\sigma^2(t):=Var\left(X(t)\right)$ such that\\
{\bf AI}: $\vf(0)=0$, $\vf(t)$ is regularly varying at $\IF$ with index $2\alpha_\IF\in(0,\eHH{2)}$
and $\vf(t)$ is twice continuously differentiable on $(0,\IF)$ with its first derivative
$\dot{\vf}(t):=\frac{{\rm d} \sigma^2}{{\rm d}t}\left(t\right)$
and second derivative
$\ddot{\vf}(t):=\frac{{\rm d^2} \sigma^2}{{\rm d}t^2}\left(t\right)$
being ultimately monotone at $\IF$.\\
{\bf AII}: $\vf(t)$ is regularly varying at $0$ with index $2\alpha_0\in(0,2]$.\\

For given $x\ge0$ and $f\in C([0,2S])$, $S>0$, let
\BQN\label{functional1}
\Gamma(x,S;f)=\sup_{t\in [0,S]}\min\left(f(t), \sup_{s\in [t+x,2S]}f(s)\right).
\EQN

Then, for given $\eta(t),t\ge0$, a centered Gaussian process with stationary increments and continuous sample paths,
we introduce
\BQN\label{gammapickands}
 \mathcal{H}_{\eta}^{\Gamma}(x,S):&=&\E{\expon{\Gamma(x,S;\sqrt{2}\eta(t)-Var(\eta(t)))}},\nonumber\\
 \mathcal{H}_{\eta}^{\Gamma}(x):&=&\lim_{S\rw\IF}\frac{\mathcal{H}_{\eta}^{\Gamma}(x,S)}{S},
\EQN
providing that the limit exists.
We note that
$\mathcal{H}_{\eta}^{\Gamma}(0)$ coincides with
the notion of {\it generalized Pickands constant},
since
\[
 \mathcal{H}_{\eta}^{\Gamma}(0,S)=\E{\expon{\sup_{t\in[0,S]}(\sqrt{2}\eta(t)-Var(\eta(t)))}};
\]
see
\cite{pit96}, \cite{DE2002}, \cite{HP2004}, \cite{DI2005} \cite{DE2014}, \cite{Harper2},
\cite{DiekerY}, \cite{DM} and \cite{Harper3}  for properties of generalized Pickands constants.
In order to simplify notation, let
\BQN\label{pickands1}
\mathcal{H}_{\eta}:=\mathcal{H}_{\eta}^{\Gamma}(0),
\quad \text{and}\quad \mathcal{H}_{\eta}([0,S]):= \mathcal{H}_{\eta}^{\Gamma}(0,S).
\EQN

Let $B_{H}(t)$ denote the standard fractional Brownian motion with mean $0$ and correlation function satisfying
$$Cov(B_H(s), B_H(t))=\frac{|s|^{2H}+|t|^{2H}-|t-s|^{2H}}{2}, \quad s,t\geq 0, H\in (0,1].$$
\COM{
In our results, the following generalized Pickands-Piterbarg type constants plays a crucial role.
 Let $E \subset \R$ be a compact set  and write $C(E)$ for the set of real-valued continuous functions defined on $E$. Let $\Gamma: C(E)\to \mathbb{R}$   be a real-valued continuous functional satisfying \\
{\bf F1:}  there exists $c>0$ such that  $\Gamma(f)\le c \sup_{{t}\in E} {f({t})}$ for any $f\in C(E)$;\\
{\bf F2:}  
$\Gamma(e_1f+e_2)=e_1\Gamma(f)+e_2$ for any $f\in C(E)$ and $e_1>0, e_2\in  \mathbb{R}$.\\
Assume that $X$ is a continuous centered Gaussian process. Define the generalized Pickands constant by
For simplicity, denote by
$$\mathcal{H}_{\eta}=\mathcal{H}_{\eta}^{\sup}, \quad \mathcal{H}_{\eta}[0,S]=\mathcal{H}_{\eta}^{\sup}([0,S]), \quad S>0.$$
One can refer to \cite{pit96}, \cite{DE2002}, \cite{HP2004}, \cite{DI2005} \cite{DE2014}, \cite{Harper2}, \cite{DiekerY}, \cite{DM} and \cite{Harper3}  for the definitions and properties of generalized Pickands constants.
 Let
$\Gamma: C([0,2S])\rw \mathbb{R}$ with $0\leq x\leq S$ and $S>0$ be a real-valued continuous functional defined by
\BQN\label{functional1}
\Gamma(x,S;f)=\sup_{t\in [0,S]}\min\left(f(t), \sup_{s\in [t+x,2S]}f(s)\right)
\EQN
for any $f\in C([0,2S])$. Apparently, $\Gamma$ satisfies {\bf F1-F2}.\\
In what follows, we display our first results.
}
Let  $ t^*=\frac{\alpha_\IF}{c(1-\alpha_\IF)}$ and $\overleftarrow{\sigma}(t), t\geq 0$
stand for the asymptotic inverse function of $\sigma$ at value of $t$.
Furthermore, let
\begin{eqnarray}\label{theta}
\Delta(u)=
\overleftarrow{\sigma}
 \left(\frac{\sqrt{2}\sigma^2(ut^*)}{u(1+ct^*)}\right).
\end{eqnarray}
In the rest of this section we tacitly assume that
\[
\varphi:=\lim_{u\rw\IF}\frac{\sigma^2(u)}{u}\in[0,\infty].
\]
\BT\label{TH4}  Let $X(t)$ be a centered Gaussian process with continuous trajectories and stationary increments satisfying
{\bf AI-AII}. 
Then for any $x\geq 0$
$$ \lim_{u\rw\IF}\pk{\frac{\mathcal{R}_\IF(u)}{\Delta(u)}>x}=\frac{\mathcal{H}_{\eta}^{\Gamma}(x)}{\mathcal{H}_{\eta}^\Gamma(0)}\in (0,1], $$
where
   \BQN\label{eta1}
\eta(t)=\left\{\begin{array}{cc}
B_{\alpha_0}(t),& \varphi=0\\
\frac{X(\varphi t)}{\sigma(\varphi)}, & \varphi\in (0,\IF)\\
B_{\alpha_\IF}(t), & \varphi=\IF.
\end{array}
\right.
 \EQN
\ET
\subsection{Finite-time horizon}
In this subsection we focus on the finite-time case,
i.e. we suppose that $T\in(0,\infty)$. 
Let $X(t), t\in [0,T]$ be a Gaussian process with stationary increments, a.s. continuous trajectories, zero-mean and variance function $\sigma^2$ satisfying \\
{\bf BI} $\sigma^2(0)=0$ and $\sigma^2\in C(0,T]$ with the first derivative being positive.\\
{\bf BII} $\sigma^2$ is regularly varying at $0$ with index $2\alpha_0\in (0,2]$.\\
Denote by
$$\Delta_1(u)=\overleftarrow{\sigma}\left(\frac{\sqrt{2}\sigma^2(T)}{u+cT}\right), \quad  \Delta_2(u)=\left(\frac{\sigma(T)}{u+cT}\right)^2,$$
and
\BQN\mathcal{P}_{B_{1/2}}^d(x)&=&\int_0^\IF e^w\pk{\sup_{t\in [x,\IF)}\min\left(\sup_{ s\in [0, t-x]}\sqrt{2}B_{1/2}(s)-\left(1+d\right)s, \sqrt{2}B_{1/2}(t)-\left(1+d\right)t\right)>w}dw\nonumber\\
&&+ \int_{-\IF}^0 e^w \pk{\sup_{t\in [x,\IF)}\left(\sqrt{2}B_{1/2}(t)-\left(1+d\right)t\right)>w}dw,
\EQN
with $B_{1/2}$ a standard Brownian motion and $d>0$.
We note that
\[
\mathcal{P}_{B_{1/2}}^d(0)= \E {\expon{\sup_{t\in [0,\infty)}(\sqrt{2}B_{1/2}(t)-(1+d)t)}},
\]
is the classical Piterbarg constant (see \cite{LongB} and reference therein) and hence
$$\mathcal{P}_{B_{1/2}}^d(0)=1+\frac{1}{d};$$
see, e.g., \cite{pit96}.
\BT\label{TH5}  Let $X(t), t\in [0,T]$ be a centered Gaussian process with continuous trajectories and stationary increments satisfying {\bf BI-BII} and $x\geq 0$.\\
i) If $t=o(\sigma^2(t)), t\to 0$, then
$$ \lim_{u\rw\IF}\pk{\frac{\mathcal{R}_T(u)}{\Delta_1(u)}>x}=\frac{\mathcal{H}_{B_{\alpha_0}}^{\Gamma}(x)}{\mathcal{H}_{B_{\alpha_0}}^\Gamma(0)}\in (0,1].$$
ii) If $\sigma^2(t)\sim at, t\to 0$, $a>0$, then
for $d=\frac{2\sigma(T)\dot{\sigma}(T)}{a}$,
$$ \lim_{u\rw\IF}\pk{\frac{\mathcal{R}_T(u)}{\Delta_1(u)}>x}=\frac{\mathcal{P}_{B_{1/2}}^{d}(x)}
{\mathcal{P}_{B_{1/2}}^{d}(0)}\in (0,1].$$
iii) If $\sigma^2(t)=o(t), t\to 0$, then
 $$ \lim_{u\rw\IF}\pk{\frac{\mathcal{R}_T(u)}{\Delta_2(u)}>x}=e^{-\frac{\dot{\sigma}(T)}{\sigma(T)}x}.$$

\ET

\section{Proofs}
In this section, we give the proof of Theorem \ref{TH4} and Theorem \ref{TH5}.
Hereafter, denote by $\mathbb{Q}$, $\mathbb{Q}_i, i=1,2,3,\dots$ positive constants that may differ from line to line and $\overline{X}:=\frac{X}{\sqrt{Var(X)}}$ for any nontrivial random variable $X$. Moreover, $f(u)\sim g(u), u\rw\IF$ means that $\lim_{u\rw\IF}\frac{f(u)}{g(u)}=1.$  In our proofs, multiple limits appear. We shall write for instance
$$b(u,S,\epsilon)\sim a(u), \quad u\rw\IF, S\rw\IF, \epsilon\rw 0 $$
to mean that
$$\lim_{\epsilon\rw 0}\lim_{S\rw\IF}\lim_{u\rw\IF}\frac{b(u,S,\epsilon)}{a(u)}=1.$$
\subsection{ Infinite-time horizon}
Observe that for any $x\geq 0$
\BQN\label{ratio1}
\pk{\frac{\mathcal{T}_\IF^*(u)-\tau_\IF^*(u)}{\Delta(u)}> x}=\frac{\pk{\mathcal{T}_\IF(u)-\tau_\IF(u)>x \Delta(u), \tau_\IF<\IF}}{\pk{\tau_\IF(u)<\IF}},
\EQN
with $\Delta(u)$ defined in (\ref{theta}).
In order to derive the limiting distribution of the above ratio, we need to derive the asymptotics of
$\pk{\mathcal{T}_\IF(u)-\tau_\IF(u)>x \Delta(u), \tau_\IF<\IF}$ and $\pk{\tau_\IF(u)<\IF}$ respectively.
Using that
$$\pk{\tau_\IF(u)<\IF}=\pk{\sup_{t\geq 0}X(t)-ct>u}=\pk{\sup_{t\geq 0}X_u(t)>m(u)},$$
where
$X_u(t)=\frac{X(ut)}{u(1+ct)}m(u)$ with
$ m(u)=\inf_{t>0}\frac{u(1+ct)}{\sigma(ut)},$
for
\BQN\label{AB}
A=\left(\frac{{\alpha_\IF}}{c(1-{\alpha_\IF})}\right)^{-\alpha_\IF}\frac{1}{1-\alpha_\IF}, \quad
B=\left(\frac{\alpha_\IF}{c(1-\alpha_\IF)}\right)^{-\alpha_\IF-2}\alpha_\IF, \quad t_u=\arg\inf_{t>0}\frac{u(1+ct)}{\sigma(ut)},
\EQN
Proposition 2 in \cite{DI2005} (or Theorems 3.1-3.3 in \cite{KrzysPeng2015}),
implies the following Lemma.

\BEL\label{TH0}
Let $X(t)$ be a centered Gaussian process with continuous trajectories and stationary increments satisfying {\bf AI-AII}. Assuming that $\lim_{u\rw\IF}\frac{\sigma^2(u)}{u}\in[0,\IF]$, we have
\BQNY
\pk{\tau_\IF(u)<\IF}\sim \mathcal{H}_{\eta}^{\Gamma}(0)\sqrt{\frac{2A\pi}{B}} \frac{u}{m(u)\Delta(u)}\Psi(m(u)),
\EQNY
where $\Delta(u)$ is defined in (\ref{theta}) and $\eta$ is defined in (\ref{eta1}).
\EEL
Thus, by (\ref{ratio1}), we are left with finding the asymptotics of
$\pk{\mathcal{T}_\IF(u)-\tau_\IF(u)>x \Delta(u), \tau_\IF<\IF}$, as $u\rw\IF$.

In the next lemma we focus on asymptotic properties of
the variance and correlation functions of relevant Gaussian processes; we refer to, e.g., \cite{KrzysPeng2015} for the proof.
Let
$$\sigma_u(t):=Var^{1/2}\left(X_u(t)\right) =Var^{1/2}\left(\frac{X(ut)}{u(1+ct)}m(u)\right), t\geq 0.$$
\BEL\label{L1}
Suppose that {\bf AI-AII} are satisfied.
 For $u$ large enough $t_u$ is unique, and  $t_u\rw t^*=\frac{\alpha_\IF}{c(1-\alpha_\IF)}$, as
 $u\to \infty$. Moreover,  for any $\delta_u>0$ with $\lim_{u\rw\IF}\delta_u=0$
\BQNY\label{V4}
\lim_{u\rw\IF}\sup_{t\in (t_u-\delta_u, t_u+\delta_u)\setminus \{t_u\}}\left|\frac{1-\sigma_u(t)}{\frac{B}{2A}(t-t_u)^2}-1\right|=0,
\EQNY
and
 \BQNY\label{cor0}
\lim_{u\rw \IF}\sup_{s\neq t, s,t \in (t_u-\delta_u, t_u+\delta_u)}\left|\frac{1-Cor\left(X(us), X(ut)\right)}{\frac{\sigma^2(u|s-t|)}{2\sigma^2(ut^*)}}-1\right|=0.
\EQNY
\EEL

\prooftheo{TH4}
Due to (\ref{ratio1}) and Lemma \ref{TH0}, we focus on the asymptotics of
\BQNY\label{decom}
&&\pk{\mathcal{T}_\IF(u)-\tau_\IF(u)>x\Delta(u), \tau_\IF<\IF}\\
 &&\quad =\pk{\exists s,t \geq 0, s-t\geq \Delta(u)x, X(t)-ct>u, X(s)-cs>u},
\EQNY
as $u\to\infty$, for any $x\geq 0$.
We have
\BQN\label{decom}
\pi_1(u)\leq \pk{\mathcal{T}_\IF(u)-\tau_\IF(u)>x\Delta(u), \tau_\IF<\IF}
\leq \pi_1(u)+\pi_2(u),
\EQN
where
\BQNY
\pi_1(u)&=&\pk{\exists s,t \in E_1(u), s-t\geq \Delta(u)x, X(t)-ct>u, X(s)-cs>u}\\
\pi_2(u)&=&\pk{\exists (s,t) \in [0,\IF)^2\setminus E_1^2(u), s-t\geq \Delta(u)x, X(t)-ct>u, X(s)-cs>u},
\EQNY
with $$E_1(u)=\left[ut_u-\frac{u\ln m(u)}{m(u)}, ut_u+\frac{u\ln m(u)}{m(u)}\right].$$
It follows that for $u>0$
\BQNY
\pi_2(u) \leq 2 \pk{\sup_{t \in [0,\IF)\setminus E_1(u)}X(t)-ct>u}.
\EQNY
Hence, following Lemma 7 in \cite{DI2005} (or Lemma 5.6 in \cite{KrzysPeng2015}), we have that
\BQN\label{neg3}
\pi_2(u) \leq 2 \pk{\sup_{t \in [0,\IF)\setminus E_1(u)}X(t)-ct>u}=o\left(\pk{\tau_\IF(u)<\IF}\right), \quad u\rw\IF.
\EQN
Thus we are left with finding the exact asymptotics of $\pi_1(u)$ as $u\rw\IF$.
Replacing $t$ by $ut_u+\Delta(u)t$ and $s$ by $ut_u+\Delta(u)s$, we rewrite
\BQN\label{pi1u}
\pi_1(u)=\pk{\exists s,t \in E_2(u), s-t\geq x, Z_u(t)>m(u), Z_u(s)>m(u)}
\EQN
with
$$ Z_u(t)=\frac{X(ut_u+\Delta(u)t)}{u(1+ct_u)+c\Delta(u)t}m(u), \quad E_2(u)=\left[-\frac{u\ln m(u)}{\Delta(u)m(u)}, \frac{u\ln m(u)}{\Delta(u)m(u)}\right].$$

Bonferroni inequality gives that for $S>x$,
\BQN\label{Bonfer5}
\Sigma_1^+(u)-\Sigma\Sigma_1(u)\leq \pi_1(u)\leq \Sigma_1^-(u) +\Sigma\Sigma_2(u),
\EQN
where
\BQNY
\Sigma_1^{\pm}(u)&=&\sum_{k=-N(u)\pm1}^{N(u)\mp1}\pk{ \exists t\in [kS, (k+1)S], s\in [t+x, (k+2)S]: \min(Z_u(t),Z_u(s))>m(u)}\\
\Sigma\Sigma_1(u)&=&\sum_{k=-N(u)-1}^{N(u)+1}\sum_{l=k+1}^{N(u)+1}\pk{  \sup_{t\in [kS, (k+1)S]} Z_u(t)>m(u),  \sup_{t\in [lS, (l+1)S]} Z_u(t)>m(u)}\\
\Sigma\Sigma_2(u)&=&\sum_{k=-N(u)-1}^{N(u)+1}\sum_{l=k+2}^{N(u)+1}\pk{  \sup_{t\in [kS, (k+1)S]} Z_u(t)>m(u),  \sup_{t\in [lS, (l+1)S]} Z_u(t)>m(u)},
\EQNY
where $N(u)=\left[\frac{u\ln m(u)}{\Delta(u)m(u)S}\right]$.
To get the asymptotics of $\pi_1(u)$, in next steps of the proof
we show that $\Sigma_1^{+}(u)\sim\Sigma_1^{-}(u)$ and
$\Sigma\Sigma_2(u)\leq \Sigma\Sigma_1(u)=o(\Sigma_1^{+}(u))$,
as $u\rw\IF$ and $S\rw\IF$.\\
{\it \underline{Asymptotics of $\Sigma_1^{\pm}(u)$}}.
Setting $$Z_{u,k}(t)=Z_u(kS+t),$$
we have that
\BQNY
\Sigma_1^{-}(u)&=&\sum_{k=-N(u)-1}^{N(u)+1}\pk{ \sup_{t\in [kS, (k+1)S]}\min\left(Z_{u}(t), \sup_{s\in [t+x, (k+2)S]}Z_{u}(s)\right)>m(u)}\\
&=&\sum_{k=-N(u)-1}^{N(u)+1}\pk{ \Gamma(x,S;Z_{u,k})>m(u)}\\
&\leq&\sum_{k=-N(u)-1}^{N(u)+1}\pk{ \Gamma(x,S;\overline{Z}_{u,k})>\frac{m(u)}{\sup_{t\in [0, 2S]}\sqrt{Var(Z_{u,k}(t))}}},
\EQNY
where
 $\Gamma$ is defined in (\ref{functional1}).
By Lemma \ref{L1}, for any $0<\epsilon<1$
$$\frac{m(u)}{\sup_{t\in [kS, (k+1)S]}\sqrt{Var(Z_{u,k}(t))}}\geq m(u)\left(1+(1-\epsilon)\frac{B}{2A}\left(|k|^*\frac{\Delta(u)}{u}S\right)^2\right):=m_{k,\epsilon}(u)$$
with $|k|^*=\min (|k|, |k+1|, |k+2|)$  as $u$ sufficiently large. Thus for $0<\epsilon<1$ and $u$ sufficiently large
\BQNY
\Sigma_1^{-}(u)\leq \sum_{k=-N(u)-1}^{N(u)+1}\pk{ \Gamma(x,S;\overline{Z}_{u,k})>m_{k,\epsilon}(u)}.
\EQNY
Using that
$$\lim_{u\rw \IF}\sup_{|k|\leq N(u)+1}\left|\frac{2\sigma^2(ut^*)}{\sigma^2(\Delta(u)) (m_{k,\epsilon}(u))^2}-1\right|=0,$$
by Lemma \ref{L1}
\BQN\label{cor2}
\lim_{u\rw \IF}\sup_{|k|\leq N(u)+1}\sup_{s\neq t, s,t \in[0,2S]}\left|m_{k,\epsilon}^2(u)\frac{1-Cor(Z_{u,k}(s), Z_{u,k}(t))}{\frac{\sigma^2(\Delta(u)|s-t|)}{\sigma^2(\Delta(u))}}-1\right|=0.
\EQN
Hence in light of  Proposition 2.3 in \cite{KEP2016},
\BQNY
\lim_{u\rw\IF}\sup_{|k|\leq N(u)+1}\left|\frac{\pk{\Gamma\left(x,S;\overline{Z}_{u,k}(t)\right)>m_{k,\epsilon}(u)}}{\Psi(m_{k,\epsilon}(u))}-\mathcal{H}_{\eta}^{\Gamma}(x,S)\right|=0,
\EQNY
where $\eta$ is defined in (\ref{eta1}).
Furthermore,
\BQN\label{pi8-}
\Sigma_1
^{-}(u)&\leq& \sum_{k=-N(u)-1}^{N(u)+1} \mathcal{H}_{\eta}^{\Gamma}(x,S)\Psi(m_{k,\epsilon}(u))\nonumber\\
&\sim&  \mathcal{H}_{\eta}^{\Gamma}(x,S)\Psi(m(u))\sum_{k=-N(u)-1}^{N(u)+1}e^{-(1-\epsilon)\frac{B}{2A} m^2(u)\left(|k|^*\frac{\Delta(u)}{u}S\right)^2}\nonumber\\
&\sim&  \frac{\mathcal{H}_{\eta}^{\Gamma}(x,S)}{S}\frac{(1-\epsilon)^{-1/2}\left(\frac{B}{2A}\right)^{-1/2}u}{m(u)\Delta(u)}\Psi(m(u))\int_{-\IF}^\IF e^{-t^2}dt\nonumber\\
&\sim&  \frac{\mathcal{H}_{\eta}^{\Gamma}(x,S)}{S}\Theta(u), \quad u\rw\IF, \epsilon\rw 0,
\EQN
where
\BQN\label{Theta}
\Theta(u):=\sqrt{\frac{2A\pi}{B}} \frac{u}{m(u)\Delta(u)}\Psi(m(u)).
\EQN
Similarly,
\BQN\label{pi8+}
\Sigma_1^{+}(u)\geq \frac{\mathcal{H}_{\eta}^{\Gamma}(x,S)}{S}\Theta(u)(1+o(1)), \quad u\rw\IF.
\EQN
{\it \underline{Upper bound for $\Sigma\Sigma_i(u), i=1,2$}}.
Similarly as in (\ref{cor2}), Lemma \ref{L1} gives that
$$\lim_{u\rw \IF}\sup_{s\neq t, s,t \in E_2(u)}\left|m^2(u)\frac{1-Cor(Z_u(s), Z_u(t))}{\frac{\sigma^2(\Delta(u)|s-t|)}{\sigma^2(\Delta(u))}}-1\right|=0,$$
By  Corollary 3.2 in \cite{KEP2016},  there exists $\mathcal{C}, \mathcal{C}_1>0$ such that for all $|k|, |l|\leq N(u)+1, l\geq k+2$,
\BQNY
\pk{\sup_{t\in [kS, (k+1)S]}\overline{Z_u}(t)>m_{k,\epsilon}(u), \sup_{t\in [lS,(l+1)S]}\overline{Z_u}(t)>m_{l,\epsilon}(u)}\leq \mathcal{C}S^2e^{-\mathcal{C}_1|k-l|^\gamma S^\gamma}\Psi(\hat{m}_{u,k,l}),
\EQNY
with $\gamma=\min(\alpha_0,\alpha_\IF)$ and $\hat{m}_{k,l}(u)=\min(m_{k,\epsilon}(u), m_{l,\epsilon}(u))$. Consequently, with aid of (\ref{pi8-}),
 \BQN\label{pi91}
 \Sigma\Sigma_2(u)&\leq &\sum_{k=-N(u)-1}^{N(u)+1}\sum_{l=k+2}^{N(u)+1}\mathcal{C}S^2e^{-\mathcal{C}_1|k-l|^\gamma S^\gamma}\Psi(\hat{m}_{k,l}(u))\nonumber\\
 &\leq&\sum_{k=-N(u)-1}^{N(u)+1}\sum_{l=k+2}^{N(u)+1}\mathcal{C}S^2e^{-\mathcal{C}_1|k-l|^\gamma S^\gamma}\left(\Psi(m_{k,\epsilon}(u))+\Psi(m_{k,\epsilon}(u))\right)\nonumber\\
 &\leq& 2\sum_{k=-N(u)-1}^{N(u)+1} \Psi(m_{k,\epsilon}(u))\sum_{|k-l|\geq 1}\mathcal{C}S^2e^{-\mathcal{C}_1|k-l|^\gamma S^\gamma}\nonumber\\
 &\leq& 2\sum_{k=-N(u)-1}^{N(u)+1} \Psi(m_{k,\epsilon}(u))\mathcal{C}S^2e^{-\mathbb{Q} S^\gamma}\nonumber\\
 &\leq&  \mathbb{Q}Se^{-\mathbb{Q} S^\gamma}\Theta(u), \quad u\rw\IF.
 \EQN
Thus, again by Proposition 2.3 in \cite{KEP2016}, taking into account
(\ref{cor2}) and noting that $\Gamma(0,S;f)=\sup_{t\in[0,S]}f(t)$, we have
 \BQNY
\lim_{u\rw\IF}\sup_{|k|\leq N(u)+1}\left|\frac{\pk{ \sup_{t\in [0,S]}\overline{Z_{u,k}}(t)>\hat{m}_{k,k+1}(u)}}{\Psi(\hat{m}_{k,k+1}(u))}-\mathcal{H}_{\eta}([0,S])\right|=0.
\EQNY
Hence,
 \BQNY
 \hat{\Sigma}(u):&=&\sum_{k=-N(u)-1}^{N(u)+1}\pk{\sup_{t\in [kS, (k+1)S]}Z_u(t)>m(u), \sup_{t\in [(k+1)S,(k+2)S]}Z_u(t)>m(u)}\\
 &=& \sum_{k=-N(u)-1}^{N(u)+1}\left(\pk{\sup_{t\in [0,S]}\overline{Z_{u,k}}(t)>\hat{m}_{k,k+1}(u)}+\pk{ \sup_{t\in [0, S]}\overline{Z_{u,k+1}}(t)>\hat{m}_{k,k+1}(u))}\right)\\
  && - \sum_{k=-N(u)-1}^{N(u)+1}\pk{\sup_{t\in [0, 2S]}\overline{Z_{u,k}}(t)>\hat{m}_{k,k+1}(u)}\\
 &\leq& \left(2\frac{\mathcal{H}_{\eta}([0,S])}{S}-\frac{\mathcal{H}_{\eta}([0,2S])}{S}\right)\Theta(u), \quad u\rw\IF,
 \EQNY
which together with (\ref{pi91}) and the fact that $$\Sigma\Sigma_1(u)=\Sigma\Sigma_2(u)+\hat{\Sigma}(u)$$ leads to
\BQN\label{sigma3}
\Sigma\Sigma_1(u)\leq \left(2\frac{\mathcal{H}_{\eta}([0,S])}{S}-\frac{\mathcal{H}_{\eta}([0,2S])}{S}+\mathbb{Q}Se^{-\mathbb{Q} S^\gamma}\right)\Theta(u), \quad u\rw\IF.
\EQN
Combination of (\ref{Bonfer5}), (\ref{pi8-}), (\ref{pi8+})-(\ref{sigma3}) yields
\BQN\label{existence}
\liminf_{u\rw\IF}\frac{\pi_1(u)}{\Theta(u)}&\geq& \frac{\mathcal{H}_{\eta}^{\Gamma}(x,S)}{S} -2\frac{\mathcal{H}_{\eta}([0,S])}{S}+\frac{\mathcal{H}_{\eta}([0,2S])}{S}-\mathbb{Q}Se^{-\mathbb{Q} S^\gamma},\nonumber\\
\limsup_{u\rw\IF}\frac{\pi_1(u)}{\Theta(u)}&\leq& \frac{\mathcal{H}_{\eta}^{\Gamma}(x,S)}{S} +\mathbb{Q}Se^{-\mathbb{Q} S^\gamma}.
\EQN
Thus under the proviso that
\BQN\label{finiteness}
\mathcal{H}_{\eta}^{\Gamma}(x)=\lim_{S\rw\IF}\frac{\mathcal{H}_{\eta}^{\Gamma}(x,S)}{S}\in (0,\IF),
\EQN letting $S\rw\IF$ in (\ref{existence}) leads to
$$\lim_{u\rw\IF}\frac{\pi_1(u)}{\Theta(u)}=\mathcal{H}_{\eta}^\Gamma(x)\in (0,\IF),$$
which combined with (\ref{ratio1})-(\ref{neg3}) establishes the claim. \\
{\it \underline{Existence of $\mathcal{H}^{\Gamma}_{\eta}(x)$}}. In order to complete the proof, we are left with proving (\ref{finiteness}).
 By (\ref{existence}), we have
$$\liminf_{S\rw\IF}\frac{\mathcal{H}_{\eta}^{\Gamma}(x,S)}{S}=\limsup_{S\rw\IF}\frac{\mathcal{H}_{\eta}^{\Gamma}(x,S)}{S}.$$
By the fact that
$$\mathcal{H}_{\eta}^{\Gamma}(x,S)\leq \mathcal{H}_{\eta}([0,S]), $$
we have
$$\limsup_{S\rw\IF}\frac{\mathcal{H}_{\eta}^{\Gamma}(x,S)}{S}\leq  \mathcal{H}_{\eta}<\IF.$$
In order to prove positivity of $\mathcal{H}_{\eta}^\Gamma(x)$, we follow the same argument as in the proof of Theorem D.2 in \cite{pit96}.
 Replacing  $\pi_1(u)$ in (\ref{pi1u})
 by
 $$\pk{\exists s\in E_2(u), t\in E_2(u)\cap \bigcup_{k\in \mathbb{Z}}[2kS, (2k+1)S], s-t\geq x, Z_u(t)>m(u), Z_u(s)>m(u)}$$
 and
 following the same arguments as for $\pi_1(u)$, we derive that for sufficiently large $\hat{S}>x$
\BQNY\label{constantl}
\liminf_{S\rw\IF}\frac{\mathcal{H}_{\eta}^{\Gamma}(x,S)}{S}\geq \frac{\mathcal{H}_{\eta}([0,\hat{S}])}{\hat{S}}-\mathbb{Q}\hat{S}e^{-\mathbb{Q} (\hat{S})^\gamma}>0,
\EQNY
where the last inequality follows by  the fact that $\mathcal{H}_{\eta}([0,S])$ is increasing with respect to $S$.
Hence,
\BQN\label{constant}
\lim_{S\rw\IF}\frac{\mathcal{H}_{\eta}^{\Gamma}(x,S)}{S}\in (0,\IF).
\EQN
This completes the proof. \QED
\subsection{Finite-time horizon}
Let $$\Delta_1(u)=\overleftarrow{\sigma}\left(\frac{\sqrt{2}\sigma^2(T)}{u+cT}\right), \quad  \Delta_2(u)=\left(\frac{\sigma(T)}{u+cT}\right)^2.$$

Observe that for any $x\geq 0$
\BQN\label{FI1}
\pk{\frac{\mathcal{T}_T^*(u)-\tau_T^*(u)}{\Delta_i(u)}> x}=\frac{\pk{\mathcal{T}_T(u)-\tau_T(u)>x \Delta_i(u), \tau_T(u)\leq T}}{\pk{\tau_T(u)\leq T}}, \quad i=1,2.
\EQN
In the following lemma
we give exact asymptotics of $\pk{\tau_T(u)\leq T}=\pk{\sup_{t\in [0,T]}X(t)-ct>u}$, referring
for the proof to Theorem 2.3
and Remark 2.4 in \cite{KEP2015} (choose $\gamma=0$).

\BEL\label{TH6}  Let $X(t), t\in [0,T]$ be a centered Gaussian process with continuous trajectories and stationary increments satisfying {\bf BI-BII}.\\
 i) If $t=o(\sigma^2(t)), t\to 0$, then
$$\pk{\tau_T(u)\leq T}\sim \mathcal{H}_{B_{\alpha_0}}^\Gamma(0)\frac{\sigma^3(T)}{\dot{\sigma}(T)}(u^2\Delta_1(u))^{-1}\Psi\left(\frac{u+cT}{\sigma(T)}\right).$$
ii) If $\sigma^2(t)\sim at, t\to 0$, $a>0$, then for $d=\frac{2\sigma(T)\dot{\sigma}(T)}{a}$,
$$\pk{\tau_T(u)\leq T}\sim \mathcal{P}_{B_{1/2}}^{d}(0)\Psi\left(\frac{u+cT}{\sigma(T)}\right).$$
iii) If $\sigma^2(t)=o(t), t\to 0$, then
$$\pk{\tau_T(u)\leq T}\sim \Psi\left(\frac{u+cT}{\sigma(T)}\right).$$
\EEL
Thus in order to derive the limit of (\ref{FI1}) as $u\rw\IF$,  it suffices to find the asymptotics of $$\pk{\mathcal{T}_T(u)-\tau_T(u)>x \Delta_i(u), \tau_T(u)\leq T}, i=1,2, \quad u\rw\IF.$$

Let
$$
\widetilde{\sigma}_u(t)=\frac{\sigma(t)}{u+ct}\frac{u+cT}{\sigma(T)}, 0\leq t\leq T.$$
\BEL\label{L3}
Suppose that {\bf BI, BII} hold. Then for $u$ sufficiently large, $\widetilde{\sigma}_u(t)$ attains its maximum over $[0,T]$ at $T$, and for any $\delta_u>0$ with $\lim_{u\rw\IF}\delta_u=0$
\BQNY
\lim_{u\rw\IF}\sup_{t\in [T-\delta_u, T]}\left|\frac{1-\widetilde{\sigma}_u(t)}{|t-T|}-\frac{\dot{\sigma}(T)}{\sigma(T)}\right|=0.
\EQNY
Moreover,
 \BQNY
\lim_{u\rw \IF}\sup_{s\neq t, s,t \in [T-\delta_u, T]}\left|\frac{1-Cor\left(X(s), X(t)\right)}{\frac{\sigma^2(|s-t|)}{2\sigma^2(T)}}-1\right|=0.
\EQNY
\EEL

\prooftheo{TH5}
By (\ref{FI1}) and Lemma \ref{TH6} it suffices to analyze asymptotics of
\BQNY
&&\pk{\mathcal{T}_T(u)-\tau_T(u)>x \Delta_j(u), \tau_T(u)\leq T}\\
&&\quad=
\pk{\mathcal{T}_T(u)-\tau_T(u)>x \Delta_j(u), \tau_T(u)\leq T, X(T)-CT<u}\\
&&\quad\quad+\pk{\mathcal{T}_T(u)-\tau_T(u)>x \Delta_j(u), \tau_T(u)\leq T,  X(T)-CT>u}\\
&&\quad=\pk{\exists 0\leq s,t\leq T, s-t\geq x\Delta_j(u),  X(t)-ct>u, X(s)-cs>u, X(T)-cT<u}\\
&&\quad \quad+ \pk{\exists 0\leq t\leq T-x\Delta_j(u),  X(t)-ct>u,  X(T)-cT>u},
\EQNY
where $j=1$ for cases i), ii) and $j=2$ for iii).

Thus, for $j=1,2$
\BQN\label{finitemain}
\sum_{i=3}^4\pi_i^{j}(u)\leq \pk{\mathcal{T}_T(u)-\tau_T(u)>x \Delta_j(u), \tau_T(u)\leq T}\leq \sum_{i=3}^4\pi_i^j(u)+\pi_5(u),
\EQN
where
\BQNY
\pi_3^{j}(u)&=&\pk{\exists T-(\ln u/u)^2\leq s,t\leq T, s-t\geq x\Delta_j(u), X(t)-ct>u, X(s)-cs>u, X(T)-cT<u},\\
\pi_4^{j}(u)&=&\pk{\exists T-(\ln u/u)^2\leq t\leq T-x\Delta_j(u),  X(t)-ct>u, X(T)-cT>u},\\
\pi_5(u)&=&\pk{\sup_{t\in [0,T-(\ln u/u)^2]}X(t)-ct>u}.
\EQNY
{\it \underline{Upper bound of $\pi_5(u)$}}. Rewrite $\pi_5(u)$ as
\BQNY
\pi_5(u)=\pk{\sup_{t\in [0,T-(\ln u/u)^2]}\frac{X(t)}{u+ct}\frac{u+cT}{\sigma(T)}>\frac{u+cT}{\sigma(T)}}
\EQNY
By the stationarity of increments of $X$ and {\bf BII}, we have
\BQNY
\mathbb{E}\left\{\left(\frac{X(t)}{u+ct}\frac{u+cT}{\sigma(T)}-\frac{X(t)}{u+ct}\frac{u+cT}{\sigma(T)}\right)^2\right\}\leq \mathbb{Q}\left(|t-s|+\sigma^2(|t-s|)\right)\leq \mathbb{Q}|t-s|^{\min(\alpha_0, 1)}
\EQNY
for $s,t\in [0,T]$.
Moreover,
following Lemma \ref{L3}, we have that for $u$ sufficiently large
\BQNY
\sup_{t\in [0,T-(\ln u/u)^2]}Var\left(\frac{X(t)}{u+ct}\frac{u+cT}{\sigma(T)}\right)\leq 1-\mathbb{Q}_1(\ln u/u)^2.
\EQNY
Consequently, by Piterbarg inequality (see, e.g., Theorem 8.1 in \cite{pit96}), for $u$ sufficiently large
\BQN\label{pi6}
\pi_5(u)\leq \mathbb{Q}_2u^{2/\min(\alpha_0, 1)} \Psi\left(\frac{u+cT}{\sigma(T)\sqrt{1-\mathbb{Q}_1(\ln u/u)^2}}\right).
\EQN
{\it \underline{Asymptotics of $\pi_i^j(u), i=3,4, j=1,2$}}. Let
$$X_u(t)=\frac{X(T-t)}{u+c(T-t)}\frac{u+cT}{\sigma(T)},\quad  0\leq t\leq T.$$
Then for $j=1,2$
\BQNY
\pi_3^j(u)&=&\pk{\exists 0\leq s,t\leq (\ln u/u)^2, t-s\geq x\Delta_j(u), X_u(t)>\frac{u+cT}{\sigma(T)}, X_u(s)>\frac{u+cT}{\sigma(T)}, X_u(0)<\frac{u+cT}{\sigma(T)}},\\
\pi_4^j(u)&=&\pk{\exists x\Delta_j(u)\leq t\leq (\ln u/u)^2, X_u(t)>\frac{u+cT}{\sigma(T)},  X_u(0)>\frac{u+cT}{\sigma(T)}}.
\EQNY
In order to derive the asymptotics of $\pi_i^j(u), i=3,4, j=1,2$, we  distinguish three scenarios: i) $t=o(\sigma^2(t))$,
ii) $\sigma^2(t)\sim at$ and iii) $\sigma^2(t)=o(t)$ as $t\rw 0$.\\
{\it $\diamond$ \underline{Case i) $t=o(\sigma^2(t))$}}. 
Clearly,
\BQN\label{pi51}
\pi_4^1(u)\leq \pk{X_u(0)>\frac{u+cT}{\sigma(T)}}=\Psi\left(\frac{u+cT}{\sigma(T)}\right).
\EQN
We are left with deriving the asymptotics of $\pi_3^1(u)$. \\
{\it\underline{Asymptotics of $\pi_3^1(u)$}}.
We note that
\BQNY
\hat{\pi}_3(u)-\pk{X_{u}(0)>\frac{u+cT}{\sigma(T)}}\leq \pi_3^1(u)\leq \hat{\pi}_3(u),
\EQNY
where
$$\hat{\pi}_3(u)=\pk{\exists 0\leq s,t\leq (\ln u/u)^2, t-s\geq x\Delta_1(u), X_u(t)>\frac{u+cT}{\sigma(T)}, X_u(s)>\frac{u+cT}{\sigma(T)}}.$$
For $S>x/2$, let
\BQN\label{Xuk}X_{u,k}(t)=X_u(\Delta_1(u)(kS+t)),\quad N_{1}(u)=\left[\frac{(\ln u)^2}{ u^2\Delta_1(u)S}\right]-1.
\EQN
Bonferroni inequality gives that
\BQNY
\Sigma_2^-(u)-\Sigma\Sigma_3(u)\leq \hat{\pi}_3(u)\leq \Sigma_2^+(u),
\EQNY
where
$$\Sigma\Sigma_3(u)=\sum_{k=0}^{N_1(u)+1}\sum_{l=k+1}^{N_1(u)+1}\pk{\sup_{t\in [kS, (k+1)S]}X_{u}(\Delta_1(u)t)>\frac{u+cT}{\sigma(T)}, \sup_{t\in [lS, (l+1)S]}X_{u}(\Delta_1(u)t)>\frac{u+cT}{\sigma(T)}},$$
and
\BQNY
\Sigma_2^{\pm}(u)&=&\sum_{k=0}^{N_1(u)\pm 1}\pk{\exists s\in [0,S], t\in [x+s, 2S], X_{u,k}(t)>\frac{u+cT}{\sigma(T)},  X_{u,k}(s)>\frac{u+cT}{\sigma(T)}}\\
&=&\sum_{k=0}^{N_1(u)\pm 1}\pk{\sup_{t\in [0,S]}\min\left( X_{u,k}(t), \sup_{s\in [x+t,2S]}X_{u,k}(s)\right)>\frac{u+cT}{\sigma(T)}}\\
&=&\sum_{k=0}^{N_1(u)\pm 1}\pk{\Gamma(x,S; X_{u,k})>\frac{u+cT}{\sigma(T)}},
\EQNY
with $\Gamma$  being defined in (\ref{functional1}).  By Lemma \ref{L3}, we have that for any $0<\epsilon<1$
$$\frac{\frac{u+cT}{\sigma(T)}}{\sup_{t\in [0,2S]}\sqrt{Var(X_{u,k}(t))}}>\frac{u+cT}{\sigma(T)}\left(1+(1-\epsilon)\frac{\dot{\sigma}(T)}{\sigma(T)}|k|\Delta_1(u)S\right):=m_{k,\epsilon,1}(u)$$

as $u$ sufficiently large. This implies that
\BQNY
\pk{\Gamma(x,S; X_{u,k})>\frac{u+cT}{\sigma(T)}}&\leq& \pk{\Gamma(x,S; \overline{X}_{u,k})>\frac{\frac{u+cT}{\sigma(T)}}{\sup_{t\in [0,2S]}\sqrt{Var(X_{u,k}(t))}}}\\
&\leq&\pk{\Gamma(x,S; \overline{X}_{u,k})>m_{k,\epsilon,1}(u)}.
\EQNY

Moreover, by Lemma \ref{L3} we have
\BQN\label{corXuk}
\lim_{u\rw \IF}\sup_{0\leq k\leq N_1(u)+1}\sup_{s\neq t, s,t \in [0,2S]}\left|m_{k,\epsilon,1}^2(u)\frac{1-Cor(X_{u,k}(t), X_{u,k}(s))}{\frac{\sigma^2(\Delta_1(u)|t-s|)}{\sigma^2(\Delta_1(u))}}-1\right|=0.
\EQN
Thus, following Proposition 2.3 in \cite{KEP2016}, we have
 \BQN\label{uniformXuk}
\sup_{|k|\leq N(u)+1}\left|\frac{\pk{\Gamma\left(x,S;\overline{X}_{u,k}(t)\right)>m_{k,\epsilon,1}(u)}}{\Psi(m_{k,\epsilon,1}(u))}-\mathcal{H}_{B_{\alpha_0}}^{\Gamma}(x,S)\right|=0,
\EQN
which combined with (\ref{constant}) implies that
\BQNY
\Sigma_2^{+}(u)&\leq& \sum_{k=0}^{N_1(u)+ 1}\mathcal{H}_{B_{\alpha_0}}^{\Gamma}(x,S)\Psi(m_{k,\epsilon,1}(u))\\
&\leq&\mathcal{H}_{B_{\alpha_0}}^{\Gamma}(x,S)\Psi\left(\frac{u+cT}{\sigma(T)}\right) \sum_{k=0}^{N_1(u)+ 1}e^{-(1-\epsilon)\frac{\dot{\sigma}(T)}{\sigma^3(T)}|k|u^2\Delta_1(u)S}\\
&\leq&\frac{\mathcal{H}_{B_{\alpha_0}}^{\Gamma}(x,S)}{S}\frac{\sigma^3(T)}{(1-\epsilon)\dot{\sigma}(T)u^2\Delta_1(u)}\Psi\left(\frac{u+cT}{\sigma(T)}\right) \int_0^\IF e^{-t}dt\\
&\sim& \mathcal{H}_{B_{\alpha_0}}^{\Gamma}(x)\Theta_1(u),\quad u\rw\IF, S\rw\IF, \epsilon\rw 0,
\EQNY
with
$$\Theta_1(u)=\frac{\sigma^3(T)}{\dot{\sigma}(T)}(u^2\Delta_1(u))^{-1}\Psi\left(\frac{u+cT}{\sigma(T)}\right).$$
Analogously,
\BQNY
\Sigma_2^{-}(u)\geq \mathcal{H}_{B_{\alpha_0}}^{\Gamma}(x)\Theta_1(u)(1+o(1)),\quad u\rw\IF, S\rw\IF.
\EQNY
 Following similar arguments as in (\ref{pi91})-(\ref{sigma3}), substituting $\eta$ by $B_{\alpha_0}$ and $\Theta(u)$ by $\Theta_1(u)$ in (\ref{pi91})-(\ref{sigma3}), we derive that
\BQNY
\Sigma\Sigma_3(u)&\leq& \left(2\frac{\mathcal{H}_{B_{\alpha_0}}([0,S])}{S}-\frac{\mathcal{H}_{B_{\alpha_0}}([0,2S])}{S}+\mathbb{Q}Se^{-\mathbb{Q} S^\gamma}\right)\Theta_1(u)\\
&=&o\left(\Theta_1(u)\right),\quad u\rw\IF, S\rw\IF.
\EQNY
Therefore,
\BQNY
\hat{\pi}_3(u)\sim \mathcal{H}_{B_{\alpha_0}}^{\Gamma}(x)\Theta_1(u),\quad u\rw\IF.
\EQNY
By the fact that
\BQNY
\hat{\pi}_3(u)-\pk{X_{u}(0)>\frac{u+cT}{\sigma(T)}}\leq \pi_3^1(u)\leq \hat{\pi}_3(u),
\EQNY
we have
$$\pi_3^1(u)\sim \mathcal{H}_{B_{\alpha_0}}^{\Gamma}(x)\Theta_1(u), \quad u\rw\IF,$$
which combined with (\ref{finitemain})-(\ref{pi51}) leads to
\BQNY
\pk{\mathcal{T}_T(u)-\tau_T(u)>x \Delta_j(u), \tau_T(u)\leq T}\sim \mathcal{H}_{B_{\alpha_0}}^{\Gamma}(x)\Theta_1(u),\quad u\rw\IF.
\EQNY
Inserting the above  and i) in Lemma \ref{TH6} to (\ref{FI1}), we establish the claim. \\
{\it $\diamond$ \underline{ii) Case $\sigma^2(t)\sim at$}}.
In this case we choose $\Delta_1(u)=\overleftarrow{\sigma}\left(\frac{\sqrt{2}\sigma^2(T)}{u+cT}\right)$ as the scaling function.\\
{\it \underline{Asymptotics of $\pi_3^1(u)$}}.
Using  notation introduced in (\ref{Xuk}), we have for $S>x$,
\BQNY
\pi_6(u)\leq \pi_3^1(u)\leq \pi_6(u)+\Sigma_3(u),
\EQNY
where
\BQNY
\pi_6(u)&=&\pk{\exists 0\leq s\leq s+x\leq t\leq S, X_{u,0}(t)>\frac{u+cT}{\sigma(T)}, X_{u,0}(s)>\frac{u+cT}{\sigma(T)}, X_{u,0}(0)<\frac{u+cT}{\sigma(T)}},\nonumber\\
&=&\pk{\Gamma'(x,S;X_{u,0})>\frac{u+cT}{\sigma(T)}, X_{u,0}(0)<\frac{u+cT}{\sigma(T)}},
\EQNY
with
\BQNY
\Gamma'(x,S;f):=\sup_{t\in [x,S]}\min\left(f(t), \sup_{s\in [0,t-x]}f(s)\right),
\EQNY
and
\BQN\label{pipi}
\Sigma_3(u)=\sum_{k=1}^{N_1(u)+1}\pk{\sup_{t\in [0,S]} X_{u,k}(t)>\frac{u+cT}{\sigma(T)}}.
\EQN
{\it \underline{Asymptotics of $\pi_6(u)$}}.
We begin with observation that
\BQNY
\pi_6(u)=\pk{\Gamma'\left(x,S;\frac{\overline{X}_{u,0}(t)}{1+(1/\sqrt{Var(X_{u,0}(t))}-1)}\right)>\frac{u+cT}{\sigma(T)}, X_{u,0}(0)<\frac{u+cT}{\sigma(T)}}.
\EQNY
We shall apply Lemma \ref{uniform} from Appendix, for which we verify assumptions {\bf D0-D2} (see Appendix).
Note that {\bf D0} holds straightforwardly. From (\ref{corXuk}) for $k=0$, we know that {\bf D1} is satisfied.
By Lemma \ref{L3}, it follows that
\BQN\label{covariance}
\lim_{u\rw\IF}\sup_{t\in [0,S]}\left|\left(\frac{u+cT}{\sigma(T)}\right)^2\left(1/\sqrt{Var(X_{u,0}(t))}-1\right)-\frac{2\sigma(T)\dot{\sigma}(T)}{a}t\right|=0.
\EQN
This implies that {\bf D2} is satisfied with $h(t)=\frac{2\sigma(T)\dot{\sigma}(T)}{a}t$. Moreover,
\BQN\label{D}
\left\{X_{u,0}(0)<\frac{u+cT}{\sigma(T)}\right\}=\left\{X_{u,0}(0)=\frac{u+cT}{\sigma(T)}-\frac{w}{\frac{u+cT}{\sigma(T)}}, w\in D\right\},
\EQN
with $D=(0,\IF)$.
Thus
\BQNY
\pi_6(u)\sim \int_0^\IF e^w \pk{\Gamma'(x,S;W)>w} dw\Psi\left(\frac{u+cT}{\sigma(T)}\right), \quad u\rw\IF.
\EQNY
with
\BQN\label{eta}
W(t)=\sqrt{2}B_{1/2}(t)-\left(1+\frac{2\sigma(T)\dot{\sigma}(T)}{a}\right)t.
\EQN
{\it \underline{Upper bound for  $\Sigma_3(u)$}}. Noting that $\Gamma(0,S; f)=\sup_{t\in [0,S]}f(t)$, by (\ref{uniformXuk}) we have
\BQN\label{pi9}
\Sigma_3(u)&\leq& \sum_{k=1}^{N_1(u)+1}\pk{\sup_{t\in [0,S]} \overline{X}_{u,k}(t)>m_{k,\epsilon,1}(u)}\nonumber\\
 &\leq&\sum_{k=1}^{N_1(u)+1}\mathcal{H}_{B_{1/2}}([0,S])\Psi\left(m_{k,\epsilon,1}(u)\right)\nonumber\\
 &\leq& \mathcal{H}_{B_{1/2}}([0,S])\Psi\left(\frac{u+cT}{\sigma(T)}\right) \sum_{k=1}^{N_1(u)+ 1}e^{-(1-\epsilon)\frac{\dot{\sigma}(T)}{\sigma^3(T)}ku^2\Delta_1(u)S}\nonumber\\
 &\leq& \mathbb{Q}  S\sum_{k=1}^{N_1(u)+ 1}e^{- \mathbb{Q}_1kS} \Psi\left(\frac{u+cT}{\sigma(T)}\right)\nonumber\\
 &\leq& \mathbb{Q}  Se^{- \mathbb{Q}_2S} \Psi\left(\frac{u+cT}{\sigma(T)}\right)=o\left(\Psi\left(\frac{u+cT}{\sigma(T)}\right)\right),\quad u\rw\IF, S\rw\IF.
\EQN
Therefore,
\BQN\label{pi41}\pi_3^1(u)\sim \int_0^\IF e^w\pk{\sup_{t\in [x,\IF)}\min\left(\sup_{ s\in [0, t-x]}\eta(s), \eta(t)\right)>w}dw \Psi\left(\frac{u+cT}{\sigma(T)}\right), \quad u\rw\IF.
\EQN
{\it\underline{ Asymptotics of  $\pi_4^1(u)$}}. Observe that
\BQN\label{pi511}
\pi_7(u)\leq \pi_4^1(u)\leq \pi_7(u)+\Sigma_3(u),
\EQN
where $\Sigma_3(u)$ is given by (\ref{pipi}) and
\BQNY
\pi_7(u)&=&\pk{\exists x\leq t\leq S, X_{u,0}(t)>\frac{u+cT}{\sigma(T)}, X_{u,0}(0)>\frac{u+cT}{\sigma(T)}}\\
&=&\pk{\sup_{ t\in [x,S] } X_{u,0}(t)>\frac{u+cT}{\sigma(T)}, X_{u,0}(0)>\frac{u+cT}{\sigma(T)}}.
\EQNY
Note that
$$ \left\{X_{u,0}(0)>\frac{u+cT}{\sigma(T)}\right\}=\left\{X_{u,0}(0)=\frac{u+cT}{\sigma(T)}-\frac{w}{\frac{u+cT}{\sigma(T)}}, w\in D\right\}
$$
with $D=(-\IF, 0)$.
 By (\ref{corXuk}), (\ref{covariance}) and (\ref{D}), applying Remark \ref{Remark} in Appendix, it follows that
\BQNY
\pi_7(u)\sim \int_{-\IF}^0 e^w \pk{\sup_{t\in [x,S]}W (t)>w}dw \Psi\left(\frac{u+cT}{\sigma(T)}\right),
\EQNY
with $W $ given in (\ref{eta}).
 Inserting the above asymptotics and (\ref{pi9}) into (\ref{pi51}) gives that
 \BQN\label{pi511}
 \pi_4^1(u)\sim \int_{-\IF}^0 e^w \pk{\sup_{t\in [x,S]}W(t)>w}dw \Psi\left(\frac{u+cT}{\sigma(T)}\right),\quad u\rw\IF.
 \EQN
 Combination of (\ref{finitemain}), (\ref{pi6}), (\ref{pi41}) and (\ref{pi511}) gives the asymptotics for $\pk{\mathcal{T}_T(u)-\tau_T(u)>x \Delta_1(u), \tau_T(u)\leq T}$, which together with ii) in Lemma \ref{TH6} establishes the claim.\\
 { \it $\diamond$ \underline{Case iii) $\sigma^2(t)=o(t)$}}.
 In this case we choose $\Delta_2(u)=\left(\frac{\sigma(T)}{u+cT}\right)^2$ as the scaling function.

{\it \underline{Asymptotics of $\pi_4^2(u)$}}. Observe that
 \BQNY
\pi_8(u)\leq \pi_4^2(u)\leq \pi_8(u)+\Sigma_4(u),
 \EQNY
 where
 \BQNY
 \pi_8(u)&=&\pk{\exists  x\leq t\leq S, \widehat{X}_{u,0}(t)>\frac{u+cT}{\sigma(T)},   \widehat{X}_{u,0}(0)>\frac{u+cT}{\sigma(T)}}\\
 &=&\pk{\sup_{t\in [x,S]}\widehat{X}_{u,0}(t)>\frac{u+cT}{\sigma(T)},  \widehat{X}_{u,0}(0)>\frac{u+cT}{\sigma(T)}}\\
 &=&\pk{\sup_{t\in [x,S]}\frac{\overline{\widehat{X}}_{u,0}(t)}{1+1/\sqrt{Var(\widehat{X}_{u,0}(t))}-1}>\frac{u+cT}{\sigma(T)},  \widehat{X}_{u,0}(0)>\frac{u+cT}{\sigma(T)}}\\
\Sigma_4(u)&=&\sum_{k=1}^{N_2(u)+1}\pk{\sup_{t\in [0,S]}\widehat{X}_{u,k}(t)>\frac{u+cT}{\sigma(T)}}\\
 &\leq& \sum_{k=1}^{N_2(u)+1}\pk{\sup_{t\in [0,S]}\overline{\widehat{X}}_{u,k}(t)>m_{k,\epsilon,2}(u)},
 \EQNY
 where
 $$\widehat{X}_{u,k}(t)=X_u(\Delta_2(u)(kS+t)),\quad N_{2}(u)=\left[\frac{(\ln u)^2}{ u^2\Delta_2(u)S}\right],$$
 $$m_{k,\epsilon,2}(u)=\frac{u+cT}{\sigma(T)}\left(1+(1-\epsilon)\frac{\dot{\sigma}(T)}{\sigma(T)}|k|\Delta_2(u)S\right).$$
 {\it \underline{Asymptotics of $\pi_8(u)$}}. In order to apply Lemma \ref{uniform} and Remark \ref{Remark} in Appendix, we check  {\bf D0-D2}. Note that  {\bf D0} hold straightforwardly.
  By Lemma \ref{L3}, it follows that
 $$\lim_{u\rw \IF}\sup_{|k|\leq N_2(u)+1}\sup_{s\neq t, s,t \in[0,S]}\left|\left(m_{k,\epsilon,2}(u)\right)^2\frac{1-Cor(\widehat{X}_{u,k}(s), \widehat{X}_{u,k}(t)}{\frac{\sigma^2(\Delta_2(u)|t-s|)}{\sigma^2(\Delta_1(u))}}-1\right|=0$$
 with
 $$\lim_{u\rw\IF}\frac{\Delta_2(u)}{\Delta_1(u)}=0.$$
 This implies that {\bf D1} holds with $\nu=0$. Lemma \ref{L3} indicates that
 \BQNY
 \lim_{u\rw\IF}\sup_{t\in [0,S]}\left|\left(\frac{u+cT}{\sigma(T)}\right)^2\left(1/\sqrt{Var(\widehat{X}_{u,0}(t))-1}\right)-\frac{\dot{\sigma}(T)}{\sigma(T)}t\right|=0.
 \EQNY
 This means that {\bf D2} holds with $h(t)=\frac{\dot{\sigma}(T)}{\sigma(T)}t$. Moreover,
 $$\left\{ \widehat{X}_{u,0}(0)>\frac{u+cT}{\sigma(T)}\right\}=\left\{ \widehat{X}_{u,0}(0)=\frac{u+cT}{\sigma(T)}-\frac{w}{ \frac{u+cT}{\sigma(T)}}, w\in D\right\},$$
 with $D=(-\IF,0)$.
Hence
 \BQNY
  \pi_8(u)&\sim&\int_{-\IF}^0 e^w \pk{\sup_{t\in [x,S]}-\frac{\dot{\sigma}(T)}{\sigma(T)}t>w}dw\Psi\left(\frac{u+cT}{\sigma(T)}\right)\\
  &=&e^{-\frac{\dot{\sigma}(T)}{\sigma(T)}x}\Psi\left(\frac{u+cT}{\sigma(T)}\right).
 \EQNY
 {\it \underline{Upper bound of $\Sigma_4(u)$}}.
Similarly, for the summands in $\Sigma_4(u)$ we can show that {\bf D0-D2} hold with $h(t)=\frac{\dot{\sigma}(T)}{\sigma(T)}t$ and $\nu=0$. Thus by Remark \ref{Remark} in Appendix, we have
\BQNY
\lim_{u\rw\IF}\sup_{|k|\leq N_2(u)+1}\left|\frac{\pk{\sup_{t\in [0,S]}\overline{\widehat{X}}_{u,k}(t)>m_{k,\epsilon,2}(u)}}{\Psi\left(m_{k,\epsilon,2}(u)\right)}-\int_{-\IF}^\IF e^w \pk{\sup_{t\in [0,S]}-\frac{\dot{\sigma}(T)}{\sigma(T)}t>w}dw\right|=0,
\EQNY
where $$\int_{-\IF}^\IF e^w \pk{\sup_{t\in [0,S]}-\frac{\dot{\sigma}(T)}{\sigma(T)}t>w}dw=1.$$
Hence
 \BQNY
 \Sigma_4(u)&\leq&\sum_{k=1}^{N_2(u)+1}\Psi\left(m_{k,\epsilon,2}(u)\right)\\
  &\leq&\Psi\left(\frac{u+cT}{\sigma(T)}\right)\sum_{k=1}^{N_2(u)+1}e^{-\mathbb{Q}|k|\Delta_2(u)u^2S}\\
  &\leq&\Psi\left(\frac{u+cT}{\sigma(T)}\right)\sum_{k=1}^{N_2(u)+1}e^{-\mathbb{Q}_1S}\\
  &\leq&e^{-\mathbb{Q}_2 S}\Psi\left(\frac{u+cT}{\sigma(T)}\right)=o\left(\Psi\left(\frac{u+cT}{\sigma(T)}\right)\right), \quad u\rw\IF, S\rw\IF.
 \EQNY
 Thus
 $$\pi_4^2(u)\sim e^{-\frac{\dot{\sigma}(T)}{\sigma(T)}x}\Psi\left(\frac{u+cT}{\sigma(T)}\right),\quad u\rw\IF.$$
{\it \underline{ Upper bound of  $\pi_3^2(u)$}}. It follows that
 \BQNY
 \pi_3^2(u)&\leq& \pk{\tau_T(u)\leq T, X_u(0)<\frac{u+cT}{\sigma(T)}}\\
 &=&  \pk{\tau_T(u)\leq T}-\pk{\tau_T(u)\leq T, X_u(0)>\frac{u+cT}{\sigma(T)}}\\
 &=& \pk{\tau_T(u)\leq T}-\pk{ X_u(0)>\frac{u+cT}{\sigma(T)}}
 \EQNY
 Applying iii) in Lemma \ref{TH6}, we have
 $$\pi_3^2(u)=o\left(\Psi\left(\frac{u+cT}{\sigma(T)}\right)\right), \quad u\rw\IF.$$
 Recalling (\ref{finitemain})-(\ref{pi6}), we conclude that
 \BQNY
 \pk{\mathcal{T}_T(u)-\tau_T(u)>x \Delta_2(u), \tau_T(u)\leq T}\sim e^{-\frac{\dot{\sigma}(T)}{\sigma(T)}x}\Psi\left(\frac{u+cT}{\sigma(T)}\right), \quad u\rw\IF,
 \EQNY
 which combined with (\ref{FI1}) and iii) in Lemma \ref{TH6} leads,  for any $x>0$, to
 $$ \lim_{u\rw\IF}\pk{\frac{\mathcal{T}_T^*(u)-\tau_T^*(u)}{\Delta_2(u)}>x}=e^{-\frac{\dot{\sigma}(T)}{\sigma(T)}x}.$$
 This completes the proof. \QED
\COM{Noting that $S>x$
$$\frac{1}{S}\expon\sup_{t\in [0,S]}\min\left(\eta(t), \sup_{s\in [t+x,2S]}\eta(s)\right)\leq \frac{1}{S}\expon\sup_{t\in [x,2S]}\eta(t)\rw 0, \quad x\rw\IF,$$
and
$$\frac{1}{S}\expon\sup_{t\in [0,S]}\min\left(\eta(t), \sup_{s\in [t+x,2S]}\eta(s)\right)\leq \frac{1}{S}\expon\sup_{t\in [0,S]}\eta(t), \quad \lim_{x\rw\IF, S>x}\mathbb{E}\left(\frac{1}{S}\expon\sup_{t\in [0,S]}\eta(t)\right)<\IF,$$
by dominated convergence theorem, we have
$$\lim_{x\rw\IF, S>x}\mathbb{E}\left(\frac{1}{S}\expon\sup_{t\in [0,S]}\min\left(\eta(t), \sup_{s\in [t+x,2S]}\eta(s)\right)\right)=0.$$
Moreover, by (\ref{pi91}) we know that $\mathbb{Q}$ in $\mathbb{Q}Se^{-\mathbb{Q} S^\gamma}$ does not depends on $x$. Thus
$$\lim_{x\rw\IF, S>x}\mathbb{Q}Se^{-\mathbb{Q} S^\gamma}=0.$$
Consequently,
$$\lim_{x\rw\IF}\mathcal{H}_{\eta}^{\Gamma}(x)=0.$$
ii) For $x=0$, in light of \cite{pit96}
\BQNY
\mathcal{P}_d^*(0)=\int_{-\IF}^\IF e^w\pk{\sup_{t\in [0,\IF)} \sqrt{2}B_{1/2}(t)-\left(1+d\right)t>w}dw=\frac{1}{2}\left(1+\sqrt{1+\frac{1}{d}}\right),
\EQNY
with $d>0$.
This implies that
$$\frac{2\mathcal{P}_d^*(0)}{1+\sqrt{1+\frac{1}{d}}}=1.$$
Moreover, using that fact that
$$\lim_{x\rw\IF} \sup_{t\in [x,\IF)}\sqrt{2}B_{1/2}(t)-\left(1+d\right)t=-\IF,$$
we have that
\BQNY
\lim_{x\rw\IF}\mathcal{P}_d^*(x)\leq\lim_{x\rw\IF}
\int_{-\IF}^\IF e^w \pk{\sup_{t\in [x,\IF)}\sqrt{2}B_{1/2}(t)-\left(1+d\right)t>w}dw=0.
\EQNY
This completes the proof. \QED}

\section*{Appendix}
In this section we give a variant of Theorem 2.1 in \cite{KEP2016}.
Let $\xi_{u,\tau_u}$
be a family of Gaussian random fields given by
 \BQN\label{xi}
 \xi_{u,\tau_u}(t)=\frac{Z_{u,\tau_u}(t)}{1+h_{u,\tau_u}(t)}, \quad t\in [0,S], \tau_u\in K_u,
 \EQN
 where $Z_{u,\tau_u}$ is a family of centered Gaussian random fields with continuous trajectories and  unit variance,
 $h_{u,\tau_u}\in C([0,S])$, $S>0$ and $K_u$ is a set of index. We investigate the asymptotics of
 $$ \pk{ \Gamma'(x,S;\xi_{u,\tau_u})>g_{u,\tau_u}, D_{u,\tau_u}}$$
 as $u\rw\IF$ where $g_{u,\tau_u}$ is a series of positive functions of $u$,
 \begin{eqnarray}
 D_{u,\tau_u}=\{Z_{u,\tau_u}(0)=g_{u,\tau_u}-\frac{w}{g_{u,\tau_u}}, w\in D\}\label{D_u}
 \end{eqnarray}
with $D=(0,\IF), (-\IF,0)$ or $D=\mathbb{R}$ and
$\Gamma': C([0,S])\rw \mathbb{R}$, $0\leq x\leq S$ is a real-valued continuous functional defined by
\BQN\label{functional11}
\Gamma'(x,S; f)=\sup_{t\in [x,S]}\min\left(f(t), \sup_{s\in [0,t-x]}f(s)\right), \quad f\in C([0,S]).
\EQN

In order to avoid trivialities, we assume that 
 $$\lim_{u\rw\IF}\pk{ \Gamma'(x,S;\xi_{u,\tau_u})>g_{u,\tau_u}, D_{u,\tau_u}}=0$$
and observe that
$\pk{ \Gamma'(x,S;\xi_{u,\tau_u})>g_{u,\tau_u}, D_{u,\tau_u}}=\pk{ \Gamma'(x,S;\xi_{u,\tau_u})>g_{u,\tau_u}}$
if $D=\R$.

As in \cite{KEP2016} (see Theorem 2.1), we impose the following assumptions:\\
 {\bf D0}: $\lim_{u\rw\IF}\inf_{\tau_u\in K_u}g_{u,\tau_u}=\IF$.\\
 {\bf D1}: There exist $\rho(t)$,  regularly varying  function at $0$ with index $2\alpha_0\in(0,2]$ and
  $b_i(u)>0, i=1,2$ satisfying $\lim_{u\rw\IF}b_i(u)=0, i=1,2$ and $ \lim_{u\rw\IF}\frac{b_1(u)}{b_2(u)}=\nu\in [0,\IF)$ such that
\BQNY
\lim_{u\rw\IF}\sup_{{\tau_u} \in K_u}\sup_{s,t\in [0,S], s\neq t}\left| \left(g_{u,\tau_u}\right)^2\frac{1-Corr(Z_{u,\tau_u}(t), Z_{u,\tau_u}(s))}{\frac{\rho(b_1(u)|t-s|)}{\rho(b_2(u))}}-1\right|=0.
\EQNY
{\bf D2}:  There exists $h\in C([0,S])$ such that
 $$\lim_{u\rw\IF}\sup_{\tau_u\in K_u}\sup_{t\in [0,S]}\left|(g_{u,\tau_u})^2h_{u,\tau_u}(t)-h(t)\right|=0.$$
 \BEL\label{uniform} Let $\xi_{u,\tau_u}$ be defined as in (\ref{xi}) and $\Gamma'$ be defined in (\ref{functional11}). Assume that {\bf D0-D2} are satisfied. Then, for $D_u,\tau_u$ defined in (\ref{D_u}) with
 $D=(0,\IF), (-\IF,0)$ or $D=\mathbb{R}$,
 $$\lim_{u\rw\IF}\sup_{\tau_u\in K_u}\left| \frac{\pk{ \Gamma'(x,S;\xi_{u,\tau_u})>g_{u,\tau_u}, D_{u,\tau_u}}}{\Psi\left(g_{u,\tau_u}\right)}-\int_{D} e^w \pk{\Gamma'(x,S;\nu^{\alpha_0} B_{\alpha_0}(t)-\nu^{2\alpha_0}|t|^{2\alpha_0}-h(t))>w}dw\right|=0.$$
 \EEL
 {\bf Proof}.  Conditioning on the event that $Z_{u,\tau_u}(0)=g_{u,\tau_u}-\frac{w}{g_{u,\tau_u}}$ and noting that
 $D_{u,\tau_u}=\{Z_{u,\tau_u}(0)=g_{u,\tau_u}-\frac{w}{g_{u,\tau_u}}, w\in D\}$,  we have
 \BQNY
 &&\pk{ \Gamma'(x,S;\xi_{u,\tau_u})>g_{u,\tau_u}, D_{u,\tau_u}}\\
 &&\quad = \frac{1}{\sqrt{2\pi}g_{u,\tau_u}}\int_{\mathbb{R}}e^{-\frac{\left(g_{u,\tau_u}-\frac{w}{g_{u,\tau_u}}\right)^2}{2}}\pk{ \Gamma'(x,S;\xi_{u,\tau_u})>g_{u,\tau_u}, D_{u,\tau_u}\Bigl|Z_{u,\tau_u}(0)=g_{u,\tau_u}-\frac{w}{g_{u,\tau_u}}}dw\\
 && \quad = \frac{e^{-\frac{(g_{u,\tau_u})^2}{2}}}{\sqrt{2\pi}g_{u,\tau_u}}\int_{D}e^{w-\frac{w^2}{2(g_{u,\tau_u})^2}}\pk{ \Gamma'(x,S;\xi_{u,\tau_u})>g_{u,\tau_u}\Bigl|Z_{u,\tau_u}(0)=g_{u,\tau_u}-\frac{w}{g_{u,\tau_u}}}dw.
 \EQNY
 Using the same procedure as in the proof of Theorem 2.1 in  \cite{KEP2016}, we can show that
 $$\int_{D}e^{w-\frac{w^2}{2(g_{u,\tau_u})^2}}\pk{ \Gamma'(x,S;\xi_{u,\tau_u})>g_{u,\tau_u}\Bigl|Z_{u,\tau_u}(0)=g_{u,\tau_u}-\frac{w}{g_{u,\tau_u}}}dw$$
 uniformly converges to
 $$\int_De^w\pk{\Gamma'(x,S;\nu^{\alpha_0} B_{\alpha_0}(t)-\nu^{2\alpha_0}|t|^{2\alpha_0}-h(t))>w}dw$$
 with respect to $\tau_u\in K_u$.
 This completes the proof. \QED

 \BRM\label{Remark}
 Lemma \ref{uniform} also holds if we substitute $\Gamma'$ by $\sup_{t\in [x,S]}f(t)$, with $f\in C[x,S]$,
 $x\geq 0$ and $D$, a measurable subset of $\mathbb{R}$ with positive Lebesgue measure.
\ERM

{\bf Acknowledgement}:
We thank Enkelejd Hashorva for discussions and comments that improved presentation of the results of this contribution.
K. D\c ebicki
was partially supported by NCN Grant No 2015/17/B/ST1/01102 (2016-2019) whereas P. Liu
was  supported by  the Swiss National Science Foundation Grant 200021-175752/1.

\bibliographystyle{plain}

 \bibliography{queue}

\end{document}